\documentclass{article}
\usepackage{amssymb}
\usepackage{amsmath}

\usepackage{cite}

\usepackage{amssymb}
\usepackage{amsfonts}
\usepackage{amscd}

\newtheorem{theorem}{Theorem}

\newtheorem{definition}[theorem]{Definition}

\newtheorem{lemma}[theorem]{Lemma}

\newtheorem{proposition}[theorem]{Proposition}

\newcommand\R{\mathbb{R}}
\newcommand\g{\frak{g}}

\newcommand\Z{\mathbb{Z}}
\newcommand\z{ \mathbb{Z}_2\times \mathbb{Z}_2}

\newcommand\h{\frak{h}}
\newcommand\m{\frak{m}}

\begin{document}

\title{\bf On Riemannian nonsymmetric spaces and flag manifolds}

\author{Abdelkader Bouyakoub \thanks{ Bouyakoub  Universit\'e d'Oran Es-S\'enia, Institut de 
Math\'ematiques, BP 1524, El Menawer 3100, ORAN, ALGERIE} \\
Michel
Goze \thanks{M.Goze@uha.fr,
{\small Universit\'{e} de Haute Alsace, F.S.T., 4, rue des Fr\`{e}res Lumi\`{e}re - 68093 MULHOUSE - FRANCE}}\\
Elisabeth Remm \thanks{corresponding author: E.Remm@uha.fr,
{\small Universit\'{e} de Haute Alsace, F.S.T., 4, rue des Fr\`{e}res Lumi\`{e}re - 68093 MULHOUSE - FRANCE.} 
This work has been supported for the three authors by AUF project MASI 2005-2006}}

\maketitle


\bigskip


\begin{abstract}
In this work we study riemannian metrics on flag manifolds adapted to the symmetries of these homogeneous nonsymmetric
spaces(. We first introduce the notion of riemannian $\Gamma $-symmetric space when $\Gamma $ is a general abelian
finite group, the symmetric case corresponding to $\Gamma =\Z_2$. We describe and study all the riemannian metrics on 
$SO(2n+1)/SO(r_1)\times SO(r_2)\times SO(r_3)\times SO(2n+1-r_1-r_2-r_3)$ for which the symmetries are isometries. 
We consider also the lorentzian case and give an example of 
a lorentzian homogeneous space which is not a symmetric space. 
\end{abstract}

\section{Introduction}
If $M$ is a homogeneous symmetric space, then at each point $x \in M$ we have a symmetry $s_x$ that is a diffeomorphim
of $M$ satisfying $s_x^2=Id$. It is equivalent to say that at every point $x \in M$ we have a subgroup
$\Gamma _x$ of {\it Diff(M)} isomorphic to $\Z_2$. 
The notion of $\Gamma$-symmetric space is a generalization of the classical notion of symmetric space by
considering a general finite abelian group of symmetries $\Gamma $ instead of $ \Z_2 $. 
The case $\Gamma =\Z_k$ was considered from the
algebraic point of view by V. Kac and the differential geometric approach was carried out by A.J. Ledger, 
M. Obata \cite{L.O} and O. Kowalski \cite{K} in terms of
$k$-symmetric spaces. A $k$-manifold is a homogeneous reductif space and the classification of 
these varieties is given by the corresponding classification of Lie algebras. The general notion of 
$\Gamma $-symmetric spaces was introduced by R. Lutz \cite{L} and was algebraically reconsidered by
  Y. Bahturin and M. Goze \cite{B.G}. In this last work the authors proved, in particular,
that a $\Gamma$-symmetric space is a homogeneous space $G/H$ and the Lie algebra $\g$ of $G$ is $\Gamma $-graded.
They give also
a classification of  $\Gamma$-symmetric spaces when $G$ is a classical simple complex Lie algebra and $\Gamma =\z$. 
We can see in particular
that the flag manifold admits such a structure. The particular case of Grassmannian manifolds comes into the framework of
symmetric manifolds. But for a general flag manifold, it is not the case. There is a great interest to study these manifolds, 
in an affine
or riemannian point of view. For example, in  loops groups theory we have to look complex algebraic homogeneous
spaces  $U_n$ and these spaces are  Grassmannians or flag manifolds. We will describe symmetries which
provide a flag manifold with a $(\z)$-symmetric structure. We then study riemannian metrics adapted to this structure, 
that is riemannian metrics for which the riemannian connection is the canonical torsion free connection 
of a homogeneous space. We have to impose in addition that the symmetries are isometries
(in the case of riemannian symmetric spaces this is a natural consequence of the very definition) 
We compute these metrics for flag manifolds and describe the associated riemannian invariants in some peculiar cases.

\section{$\Gamma $-symmetric spaces}
   In this section we recall some basical notions (see
\cite{B.G} for more details).
\subsection{Definition}
Let $\Gamma $ be a finite abelian group. A  $\Gamma $-symmetric space is a triple $(G,H,\Gamma_G )$ where
$G$ is a connected Lie group, $H$ a closed subgroup of $G$ and $\Gamma_G $ an abelian finite subgroup of the
group of automorphisms of $G$ isomorphic to $\Gamma $:
$$\Gamma_G =\{\rho _\gamma \in Aut(G), \ \gamma \in \Gamma \}$$
such that $H$ lies between $G_\Gamma $ the closed subgroup of $G$ consisting of all elements left fixed by the 
automorphisms of $\Gamma_G $ and the identity component of $G_\Gamma $.
The elements of $\Gamma_G $ satisfy :
$$
\left\{
\begin{array}{l}
\rho _{\gamma _1}\circ \rho _{\gamma _2}=\rho _{\gamma _1\gamma _2}, \\
\\
 \rho _e=Id  \ \ {\mbox{\rm where}} \ $e$ \ {\mbox{\rm is the unit element of }} \ G, \\
\\
\rho _\gamma (g)=g \ ,  \forall \gamma \in \Gamma  \Longleftrightarrow g\in H.
\end{array}
\right.
$$
We also suppose that
$H$ does not contain any proper normal subgroup of $G$. 

\subsection{$\Gamma $-symmetries on the homogeneous space $M=G/H$}

Given a $\Gamma $-symmetric space $(G,H,\Gamma_G )$ we construct for each point $x$ of $M=G/H$ a subgroup
$\Gamma _x$ of ${\mbox {\it Diff}}(M)$, the group of  diffeomorphisms of $M,$ isomorphic to $\Gamma $ which has  $x$ as an isoled fixed point. 
We denote by $\bar g$ the class
of $g \in G$ in $M$ and $e$  the identity of $G$.
We consider
$$\Gamma _{\bar e}=\{s  _{(\gamma ,\bar e)} \in {\mbox {\it Diff}}(M), \ \gamma \in \Gamma \}$$
with $s  _{(\gamma ,\bar e)}(\bar g)=\overline{\rho _\gamma (g)}.$

\medskip

\noindent In another point $x=\overline{g_0}$ of $M$ we have
$$\Gamma _{x}=\{s  _{(\gamma ,x)} \in {\mbox {\it Diff}}(M), \ \gamma \in \Gamma \}$$
with $s  _{(\gamma ,\bar {g_0})}(y)=g_0(s  _{(\gamma ,\bar e)})(g_0^{-1}y)$.
All these subgroups $ \Gamma _{x}$ of ${\mbox {\it Diff}}(M)$ are isomorphic to $\Gamma $.

\medskip

\noindent Since for every $x \in M$ and $\gamma \in \Gamma $, the map $ s _ {(\gamma ,x)} $ is a diffeomorphism of $M,$ 
such that $s_ {(\gamma ,x)}(x)=x$ the
tangent linear map $(Ts_ {(\gamma ,x)})_x$ is in $GL(T_xM)$. 
Thus, for every $x \in M$, we obtain a linear representation
$$S _x:\Gamma \longrightarrow GL(T_xM)$$
defined by
$$S _x(\gamma )=(T s _ {(\gamma ,x)})_x$$
and $S (\gamma )$ can be considered as a $(1,1)$-type tensor on $M$ satisfying 

1.  For every $\gamma \in \Gamma $, the map $x \in M \longrightarrow S _x (\gamma )$ is of
class $\cal{C}^{\infty },$

2. For every $x \in M$, $\{X_x \in T_x(M)\ {\mbox{\rm such \ that}} \ S  _{x}(\gamma )(X_x)=X_x, \ 
\forall \gamma \}=\{0\}$.

\noindent If we denote by
$$\check \Gamma _x=\{S _x(\gamma ), \gamma \in \Gamma \}$$
then $\check \Gamma _x$ is a subgroup of $GL(n,T_x(M))$ isomorphic to $\Gamma $.

\subsection{$\Gamma $-grading of the Lie algebra of $G$}
Let $\g$ be the Lie algebra of $G$. Each automorphism $\rho _\gamma $ of $G$ induces an automorphism $\tau _\gamma $
of $\g$. Let $\check  \Gamma $ be the set of all these automorphisms $\tau _\gamma $. 
Then $\check  \Gamma $ is a finite abelian
subgroup of $Aut(\g)$ isomorphic to $\Gamma$ and $\g$ is graded by $\Gamma $ 
that is
$$\g=\oplus _{\gamma \in \Gamma } \ \g_\gamma $$
with $[\g_{\gamma _1},\g_{\gamma _2}] \subset \g_{\gamma _1\gamma _2}.$ The group $\check  \Gamma $ is 
canonically isomorphic to
the dual group of $\Gamma .$  
Conversely every $\Gamma $-grading of $\g$ defines a $\Gamma $-symmetric space $(G,H,\Gamma_G )$ where
$G$ is a Lie group corresponding to $\g$ and the Lie algebra of $H$ is the component $\g_e$ corresponding
to the identity of $\Gamma $.

\subsection{Canonical connections of a $\Gamma $-symmetric space}

If $(G,H,\Gamma_G )$ is a $\Gamma $-symmetric space, the homogeneous space $M=G/H$ is reductive. In fact the Lie algebra
$\g$ being $\Gamma $-graded we have $\g=\oplus \ \g_\gamma= \g_e \oplus  \frak{m}$ with 
$\frak{m}=\oplus _{\gamma \in \Gamma , \gamma \neq e} \ \g_\gamma $ and  $[\g_e,\m] \subset \m$. 
If we suppose  $H$  connected, this last relation
means that $ad(H)(\m )\subset \m$. If $ad(H)(\m )=\m$, then any
connection on $G/H$ invariant by left translations of $G$ is defined by the $\g_e$-component $\omega $ of the
canonical $1$-form  $\theta $ of $G$. In this case the curvature $\Omega $ is given by
$$\Omega (X,Y)=-\frac{1}{2}[X,Y]_{\g_e}$$
for every $X,Y \in \m$. Moreover the Lie algebra of the holonomy group in $\bar e$ is generated by all elements of the form
$[X,Y]_{\g_e}$, $X,Y \in \m.$ This connection is called \cite{K.N} the canonical connection of the principal 
fibered bundle $G(G/H,H)$. Its torsion and 
curvature are given at the origin $\bar e$ of $G/H$ by
$$
\begin{array}{l}
T(X,Y)_{\bar e}=-[X,Y]_{\m}  \\
\\
R(X,Y)_{\bar e}=-[X,Y]_{g_e}
\end{array}
$$
for all $X,Y \in \m.$

\noindent If $\Gamma =\Z/2\Z$, that is if $(G,H,\Gamma_G )$ is a symmetric space, then the canonical
connection $\nabla$ on $M=G/H$ is  torsion free. In all the other cases, for example when $\Gamma $ is the Klein group, 
the torsion $T$ of $\nabla$ does not vanish. We consider then the connection $\overline \nabla$ given by
$$ \overline \nabla = \nabla -T.$$
This connection is  torsion free. Its curvature tensor writes
$$
\begin{array}{ll}
(R_{\overline \nabla}(X,Y)(Z))_{\bar e}= & \frac{1}{4}[X,[Y,Z]_{\m}]_{\m} -\frac{1}{4}[Y,[X,Z]_{\m}]_{\m} 
 -\frac{1}{2} [[X,Y]_{\m},Z]_{\m} \\ 
& \\
& -[[X,Y]_{\g_e},Z]_{\m} 
\end{array}
 $$
for all $X,Y,Z \in \m$ while the curvature of $\nabla$ is given by 
$$
(R_{ \nabla}(X,Y)(Z))_{\bar e}= -[[X,Y]_{\g_e},Z]_{\m}. $$
The geodesics of $\nabla$ and $\overline \nabla$ are the same. The connection 
$\overline  \nabla$ is called
the torsionfree canonical connection.
We can note that the canonical connection satisfies also
$$
\begin{array}{l}
\nabla T=0\\
\nabla R_{\nabla}=0.
\end{array}
$$
Moreover the symmetries $S_x(\gamma )$ are affine transformations with respect to $\nabla$.

\section{Riemannian $\Gamma$-symmetric spaces}

\subsection{Riemannian symmetric space}

Let $M=G/H$ a homogeneous symmetric space, where $G$ is a connected Lie group. We denote by $0$ the coset $H$ of $M$, 
that is the class on $G/H$ of the identity $1$ of $G$. The Lie algebra $\g$ of $G$ is $\Z_2$-graded
$$\g=\g_e \oplus \g_a$$
where $ \Z _2 =\left\{ e,a \right\}$ and this decomposition is $ad(H)$-invariant. The Lie algebra of $H$ is $\g_e$ and the 
tangent space at $0$ $T_0 M$ is identified to $\g_a$.

Every $G$-invariant metric $g$ on $G/H$ is given by an $ad(H)$-invariant non degenerate symmetric bilinear form
$B$ on $\g/\g_e$ by
$$B_\g(\bar{X},\bar{Y})=g(X,Y)$$ 
for $X,Y \in \g$ and $\bar{X}$ the class of $X$ in $\g / \g_e$. We identify $X\in \g$ with the projection on $M$
of the associated left invariant vector field on $G$. Moreover $g$ is a riemannian metric if and only if $B$
is positive definite. 
The identification of $\g/ \g_e$ with $\g_a$ permits to consider $B$ as  a non degenerate bilinear form
on $\g_a$. This form satisfies
$$B(X,[Y,Z]_{\g_a})=B(X,0)$$
for all $Y,Z \in \g_a$ because $[\g_a,\g_a] \subset \g_e.$
Then $B(X,[Y,Z]_{\g_a})+B([Y,X],Z )=0$ for all $X,Y,Z \in \m=\g_a$ and $M=(G/H,g)$ is naturally reductive. This 
means that the riemannnian connection of $G$ coincides with the canonical torsion free connection $\overline{\nabla }$
of $M$ and the symmetries $S_x \in \Gamma_x$ for all $x \in M$ are isometric. Conversely let $g$ be a metric on $G/H$
such that for each $x \in M$ $S_x$ is an isometry. If $ad(H)$ is a compact subgroup of $GL(\g)$, then there exists 
an $ad(H)$-invariant inner product $\tilde{B}$ on $\g$ such that 
$$\begin{array}{l}
1)\tilde{B}(\g_e,\g_a)=0 \\
2) \tilde{B}\mid_{\g_a}=B \mbox{ {\rm induces  the  riemannian metric $g$ on $G/H$}}
\end{array}
$$
Since $[\g_a,\g_a]\subset \g_e$ the naturally reductivity is obvious and the riemannian connection coincides with 
$\overline{\nabla }$.

Recall that if $G$ is a semi-simple Lie group then $B$ is neither but the restriction  to $\g_a$ of the Killing-Cartan 
form $\tilde{B}$ on $G$ that is 
$$B(X,Y)=tr(adX \circ adY)$$
for all $X,Y \in \m$.

\subsection{Riemannian $\Gamma $-symmetric spaces}

Let $\Gamma $ be a finite abelian group not isomorphic to $\Z_2$ and $g$ any $G$-invariant metric on a $\Gamma $-symmetric
space $M=G/H$.
Let us quppose that the symmetries $S_x$ are isometries for $g$. As $\Gamma $ is not isomorphic to $\z$, this property  doesn't imply in general the coincidence of the 
associated Levi-Civita connection and $\overline{\nabla }$. 

\begin{definition}

Let $(G,H,\Gamma_G)$ be a $\Gamma$-symmetric space and $g$ a $G$-invariant metric on $M$. 
We say that $(M,g)$ is a riemannian $\Gamma $-symmetric space if the symmetries $S_x$ are isometries for all $x \in M.$
\end{definition}

\begin{lemma}
Let $(G,H,\Gamma_G)$ a $\Gamma $-symmetric space and $\g=\oplus_{\gamma \in \Gamma }\g_{\gamma}$ the associated 
$\Gamma $-grading of the Lie algebra $\g$ of $G$. Then for every $\gamma \in \Gamma$
$$ad(H)\g_{\gamma} \subset \g_{\gamma}$$
\end{lemma}
{\it Proof. }Let $X$ be in $\g_\gamma $. For every  $\tau _\alpha   \in  \check \Gamma $, we have
$\tau _\alpha  (X)=\lambda (\gamma,\alpha  ) X$ with $\lambda (\gamma ,\alpha )=\pm 1.$ Then
$$\tau _\alpha  (ad(h)(X))=ad(\rho _\alpha (h))(\tau _\alpha  (X))=\lambda(\gamma,\alpha  ) ad(h)(X)$$
because all the elements of $H$ are invariant by the automorphisms $\rho _\alpha $. This proves that $ad(h)X \in \g_\gamma .$

\begin{proposition}
If $ad(H)$ is a compact subgroup of $GL(\g)$ and $g$ a $G$-invariant metric on the $\Gamma $-symmetric space
$M=G/H$ then there exits an $ad(H)$-inner product $\tilde{B}$ on $\g$ such that 

$$\begin{array}{l}
1)\tilde{B}(\g_{\gamma},\g_{\gamma'})=0 
\mbox{ {\rm for }} \gamma \neq \gamma' {\mbox{ \rm in}} \  \Gamma   \\
2) \tilde{B}\mid_{\g_a}=B \mbox{ {\rm induces  the riemannian metric $g$ on $G/H$}}
\end{array}
$$
\end{proposition}

\noindent {\it Proof}. Since each homogeneous component $\g_{\gamma}$ is invariant by $ad(H)$, there exists
an inner product $B$ on $\g$ which is $ad(\g_e)$-invariant and which defines $g$. As the symmetries $s_(\gamma,x)$
are isometries, we deduce that the automorphisms $\tau_\gamma$ are isometries for $\tilde{B}$. 
If $X \in \g_{\gamma},Y \in \g_{{\gamma}'}$, there exits $\alpha \in \Gamma$ such that 
$$\tau_{\alpha}(X)=\lambda (\alpha ,\gamma )X, \ \tau_{\alpha}(Y)=\lambda (\alpha ,\gamma' )Y$$
with $\lambda (\alpha ,\gamma )\lambda (\alpha ,\gamma' )=-1.$ Thus
$$\tilde{B}(X,Y)=\tilde{B}(\tau_\alpha (X),\tau_\alpha (Y))=-\tilde{B}(X,Y) \mbox{\rm{and }} \tilde{B}(X,Y)=0.$$

\medskip

\noindent {\bf Example.} Let us consider the $\z$-symmetric space $$(S0(5);SO(2) \times SO(2) \times SO(1), \Gamma_G)$$  
where $\Gamma_G$ is defined as follows.
One writes a general element of $so(5)$ by
$$so(5)=\left\{
\left(
\begin{array}{ccccc}
0 & x_1 & a_1 & a_2 & b_1 \\
-x_1 & 0 & a_3 & a_4 & b_2 \\
-a_1 & -a_3 & 0 & x_2 & c_1 \\
-a_2 & -a_4 & -x_2 & 0 & c_2 \\
-b_1 & -b_2 & -c_1 & -c_2 & 0
\end{array}
\right),x_i,a_i,b_i,c_i \in \R\right\}.
$$
We put 
$$\begin{array}{l}
\g_e=\left\{ X \in so(5)\,  /  \, a_i=b_i=c_i=0 \right\}, \\
\g_a=\left\{ X \in so(5)\,  /  \, x_i=b_i=c_i=0 \right\}, \\
\g_b=\left\{ X \in so(5)\,  /  \, x_i=b_i=c_i=0 \right\}, \\
\g_c=\left\{ X \in so(5)\,  /  \, x_i=a_i=b_i=0 \right\}.
\end{array}$$
If $\z=\left\{ e,a,b,c \right\}$, then $so(5)=\g_e \oplus \g_a \oplus g_b \oplus g_c$ is a $(\z)$-grading.
In this case
$$\check{\Gamma}=\left\{ \tau_e,\tau_a,  \tau_b,  \tau_c  \right\}$$
with $\tau_e=id, \tau_a(X)=X$ for $X \in \g_e \oplus \g_a$, $\tau_a(X)=-X$ for $X \in \g_b \oplus \g_c$, 
$\tau_b(X)=X$ for $X \in \g_e \oplus \g_b$, $\tau_b(X)=-X$ for $X \in \g_a \oplus \g_c$ and 
$\tau_c(X)=X$ for $X \in \g_e \oplus \g_c$, $\tau_c(X)=-X$ for $X \in \g_a \oplus \g_b$.
Since $G=SO(5)$ is connected, this grading gives a $(\z)$-symmetric structure on 
$M= S0(5)/SO(2) \times SO(2) \times SO(1)$ and $\g_e$ is the Lie algebra of $H=SO(2) \times SO(2) \times SO(1)$.
We denote by $\left\{ \left\{ X_1,X_2 \right\}, \left\{ A_1,A_2,A_3,A_4 \right\}, 
\left\{ B_1,B_2 \right\}, \left\{ C_1,C_2 \right\} \right\}$ the basis of $so(5)$ where each big letter corresponds to 
the matrix of $so(5)$ with the small letter equal to $1$ and other  coefficients are zero. This basis is 
adapted to the grading. Let us denote by 
$\left\{ \omega_1, \omega_2, \alpha_1, \alpha_2, \alpha_3, \alpha_4, \beta_1, \beta_2, \gamma_1, \gamma_2  \right\}$  
the dual basis.
 
\medskip

\noindent {\it
Every $ad(H)$-invariant symmetric bilinear form $B$ on $\m=\g_a \oplus \g_b \oplus \g_c$ such that 
$B(\g_\gamma,\g_{\gamma'})=0$ for $\gamma \neq \gamma'$ in $\Gamma$ is written
$$B=t(\alpha_1^2+ \alpha_2^2+\alpha_3^2+\alpha_4^2)+u(\alpha_1 \alpha_4-\alpha_2 \alpha_3 )+
v(\beta_1^2 +\beta_2^2)+w( \gamma_1^2+ \gamma_2^2).$$  
}

\noindent In fact, since $H$ is connected the bilinear product $B$ is $ad(H)$-invariant if and only if 
$$B([X,Y],Z)+B(Y,[X,Z])=0$$
for $Y,Z \in \m$ and $X \in \h=\g_e$.

The brackets of $so(5)$ with respect to the basis $\left\{ X_i,A_i,B_i,C_i \right\}$ are summarized
in the follwing table
$$\begin{array}{r|r|r|r|r|r|r|r|r|r|r}
    & X_1 & X_2 & A_1  &  A_2 & A_3 & A_4 & B_1 & B_2 & C_1 & C_2 \\
\hline
X_1 & 0   & 0   & -A_3 & -A_4 & A_1 & A_2 &-B_2 & B_1 & 0   & 0   \\
X_2 &     & 0   & -A_2 &  A_1 &-A_4 & A_3 &   0 &   0 &-C_2 & C_1 \\
A_1 &     &     & 0    & -X_2 &-X_1 & 0   &-C_1 &  0  & B_1 & 0   \\
A_2 &     &     &      & 0    &0    & -X_1&-C_2 &  0  & 0   & B_1 \\
A_3 &     &     &      &      & 0   & -X_2& 0   &-C_1 & B_2 & 0   \\
A_3 &     &     &      &      &    &  0   & 0   &-C_2 & 0   & B_2 \\
B_1 &     &     &      &      &    &      & 0   &-X_1 &-A_1 &-A_2 \\
B_2 &     &     &      &      &    &      &     & 0   &-A_3 &-A_4 \\
C_1 &     &     &      &      &    &      &     &     & 0   &-X_2 \\
C_2 &     &     &      &      &    &      &     &     &     &  0  \\
\end{array}$$
The identity $B([X_i,A_j],A_j)=0$ implies 
$$ B(A_1,A_3)=B(A_1,A_2)=B(A_2,A_4)=B(A_3,A_4)=0,$$
$$B(B_1,B_2)=B(C_1,C_2)=0.$$
The identity $B([X_2,A_i],A_j)+B(A_i,[X_2,A_j])=0$ gives for $i \neq j$
$$
\begin{array}{l}
B(A_2,A_3)+B(A_1,A_4)=0, \\
-B(A_3,A_3)+B(A_1,A_1)=0 \\
-B(A_4,A_4)+B(A_2,A_2)=0 \\
-B(A_2,A_2)+B(A_1,A_1)=0 \\
\end{array}
$$ 
In the same way we find
$$
\begin{array}{l}
B(B_1,B_1)=B(B_2,B_2), \\
B(C_1,C_1)=B(C_2,C_2)
\end{array}
$$ 
this gives 
$$B=t(\alpha_1^2+ \alpha_2^2+\alpha_3^2+\alpha_4^2)+u(\alpha_1 \alpha_4-\alpha_2 \alpha_3 )+
v(\beta_1^2 +\beta_2^2)+w( \gamma_1^2+ \gamma_2^2).$$  

\medskip

\noindent {\it 
The metric $g$ on $$S0(5)/SO(2) \times SO(2) \times SO(1)$$  associated to $B$ is naturally reductive if and only if $t=v=w$ and $u=0$.
}

\noindent In fact, if $\g$ is naturally reductive then $B$ satisfies 
$$B(X,[Z,Y]_m)+B([Z,X]_\m,Y)=0$$ 
for every $X,Y,Z \in m.$ In particular $B(A_1,[B_2,C_2]_\m)+B([C_2,A_1]_\m,B_2)=0$ gives 
$-B(A_1,A_4)+B(0,B_2)=0$ and $u=0.$
Similarly $B(A_1,[B_1,C_1])+B([B_1,A_1],C_1)=0$ gives 
$-B(A_1,A_1)+B(C_1,C_1)=0$ that is $t=w,$
and $B(B_1,[A_1,C_1])+B([A_1,B_1],C_1)=0$ gives 
$B(B_1,B_1)-B(C_1,C_1)=0$ that is $v=w.$

\begin{proposition}
The riemannian connection $\nabla_g$ of the metric $g$ on $S0(5)/SO(2) \times SO(2) \times SO(1)$  coincides with 
the canonical torsion free connection $\overline{\nabla }$ if and only if 
$B=\sum_{i=1}^4 \alpha_i^2 + \sum_{i=1}^2 \beta_i^2 +\sum_{i=1}^2 \gamma_i^2.$
\end{proposition}

\noindent {\it Remark} If $g$ is a $G$-invariant metric on $G/H$ such that its connection $\nabla_g$
is equal to $\overline{\nabla }$ the bilinear form $B$ is naturally reductive. In the previous example,
since $G$ is a simple Lie group, this inner product $B$ is the restriction to $\m$ of the Kiling-Cartan
form $K$ of $G$.
$$B(X,Y)=K(X,Y)=tr(adX \circ ad Y).$$
Then the homogeneous component $\g_\gamma$ are pairwise orthogonal and the $\tau_\gamma$ are isometries. But 
it is not the case in general.

\bigskip

Let us return to the general case.

\begin{definition}
Let $(G,H,\Gamma_G,g)$ a riemannian $\Gamma$-symmetric space. We say that $g$ is adapted to the 
$\Gamma $-structure if the Levi-Civita connection coincides with the canonical one.
\end{definition}

\begin{proposition}
Every riemannian $\Gamma$-symmetric space with adapted riemannian connection is naturally reductive with 
respect to the decomposition $\g=\g_e \oplus \m$ with $\m= \oplus_{\Gamma \neq e}\g_\gamma.$
\end{proposition}
{\it Proof.} Any $G$-invariant riemannian metric $g$ on a reductive homogeneous space $G/H$ 
with an $ad(H)$-invariant decomposition $\g=\g_e \oplus \m$ corresponds to an $ad(H)$-invariant 
non degenerate symmetric bilinear form $B_{\m}$ on $\m$. Since $M=G/H$ is a riemannian 
$\Gamma$-symmetric space,  its $G$-invariant riemannian metric $g$ is parallel with respect to the canonical
torsionless connection $\overline \nabla.$ Then from \cite{K.N} Theorem 3.3 the riemannian connnection
of $g$ and $\overline \nabla$  coincides on $G/H$ if and only if 
$B_{\m}$ satisfies
$$ B_{\m}(X,[Y,Z]_{m})+B_{\m}([Y,Z]_{\m},X)=0$$
for all $X,Y,Z \in \m$. This means that $(G/H,g)$ is naturally reductive. 

\subsection{Irreducible riemannian $\Gamma $-symmetric spaces}

Let $(G,H,\Gamma_G)$ a $\Gamma $-symmetric space. Since $G/H$ is a reductible homogeneous space with an
$ad\, H$ invariant decomposition $\g=\g_e \oplus \m$ then the Lie algebra of the holonomy group of 
$\nabla$ is spanned by the endomorphisms of $\m$ given by $R(X,Y)_0$ for all $X,Y \in \m$. Recall that
$(R(X,Y)Z)_0=-[[X,Y]_{\h},Z]$ for all $X,Y,Z \in \m$. In particular we have 
$R(X,Y)_0=0 $ as soon as $X \in \g_\gamma , Y \in \_{\gamma '} $ with $\gamma ,\gamma' \neq e$. 
For example if $\Gamma=\z $ then $\g=\g_e \oplus \g_a \oplus \g_b \oplus \g_c$ and 
$R(\g_a,\g_b)_0=R(\g_a,\g_c)_0=R(\g_b,\g_c)_0=0$. 

\begin{lemma}
Let $\g$ is a simple Lie algebra $\z$-graded. Then 
$$[\g_a,\g_a] \oplus [\g_b ,\g_b] \oplus [\g_c ,\g_c]=\g_e.$$
\end{lemma}
{\it Proof.} Let $U$ denote $[\g_a,\g_a] \oplus [\g_b ,\g_b] \oplus [\g_c ,\g_c]$. 
Then $I=U \oplus \g_a\oplus \g_b \oplus \g_c$ is an ideal of $\g$. 
In fact $X \in I$ is decomposed as $X_{U}+X_a +X_b +X_c$. The main point is to prove that
$[X_{U}, Y]$ is in $I$ for any $Y \in \g_e$. But $X_U$ is decomposed as $ [X_a,Y_a]+ [X_b,Y_b]+[X_c,Y_c]$.
The Jacobi identity shows that $[[X_a, Y_a],Y] \in [\g_a,\g_a].$
It is similary for the other components. Then $I$ is an ideal of $\g$ which is simple so $U=\g_e$.

\medskip

Remark that in any case , as soon as $\Gamma$ is not $\Z_2$ the representation $ad \, \g_e$ is not ireeducible on $\m$.
In fact  each component $\g_\gamma $ is an invariant subspace of $\m$. 

\begin{definition}
The representation $ad \, \g_e$ on $\m$ is called $\Gamma$-irreducible if $\m$ can not be written
$\m=\m_1 \oplus \m_2$ with $\g_e \oplus m_1$ and $\g_e \oplus m_2$ are $\Gamma $-graded Lie algebras.
\end{definition}
{\bf Example.} Let $\g_1$ be a simple Lie algebra and $\g=\g_1 \oplus \g_1 \oplus \g_1 \oplus \g_1.$
Let $\sigma_1, \sigma_2, \sigma_3$ the automorphisms of $\g$ given by
$$
\left\{
\begin{array}{l}
\sigma_1(X_1,X_2,X_3,X_4)= (X_2,X_1,X_3,X_4), \\
\sigma_2 (X_1,X_2,X_3,X_4)= (X_1,X_2,X_4,X_3), \\
\sigma_3=\sigma_1 \circ \sigma_2.
\end{array}
\right.
$$
They define a $(\z)$-graduation on $\g$ and we have
$\g_e=\left\{ (X,X,Y,Y) \right\},\g_a=\left\{ (0,0,Y,-Y) \right\},\g_b=\left\{ (X,-X,0,0) \right\}$ and 
$\g_c=\left\{ (0,0,0,0) \right\}$ with $X,Y \in \g_1$. In particular $\g_a$ is isomorphic to $\g_1$ 
so we have $[\g_e , \g_a]=\g_a$ and since $\g_1$ is simple we can not have $\g_a=\g_a^1 +\g_a^2$ with 
$[\g_e , \g_a^i]=\g_a^i$ for $i=1,2$. Then $\g$ is $(\z)$-graded and this decomposition is $(\z)$-irreducible. 

\medskip

Suppose now that $\g$ is a simple Lie algebra. Let $K$ be the Killing-Cartan form of $\g$. It is invariant by 
all automorphisms of $\g$. In particular
$$K(\tau_\gamma X,\tau_\gamma Y)=K(X,Y)$$
for any $\tau_\gamma  \in \check \Gamma.$ If $X \in \g_\alpha$ and $Y \in \g_\beta, \, \alpha \neq \beta $ there
exists $\gamma \in \Gamma$ such that $\tau_\gamma X=\lambda(\alpha,\gamma)X$ 
and $\tau_\gamma Y=\lambda(\beta,\gamma)Y$ with $\lambda(\alpha ,\gamma ) \lambda(\beta ,\gamma )\neq 1$. Thus $K(X,Y)=0$ and the homogeneous
components $\g_\gamma $ are pairewise orthogonal with respect to $K$. Moreover $K_\gamma =K \! \! \mid _{\g_\gamma }$
is a nondegenerate bilinear form. Since $\g$ is a simple Lie algebra, there exists an $ad \, \g_e$-invariant inner product 
$\tilde{B}$ on $\g$ such that the restriction $B=\tilde{B} \! \! \mid _\m$ to $\m$ defines a riemannian $\Gamma$-symmetric 
structure on $G/H$. This means that $\tilde{B}(\g_\gamma ,\g_{\gamma'})=0 $ for $\gamma \neq \gamma ' \in \Gamma.$
We consider an orthogonal basis of $\tilde{B}.$ For each $X \in \g_e,$ $ad \, X$ is expressed by a 
skew-symmetric matrix $(a_{ij}(X))$ and $K(X,X)=\sum_{i,j}a_{ij}(X)a_{ji}(X)<0$. So $K$ is negative-definite on $\g_e$.

Let $K_\gamma $ and $B_\gamma $ be the restrictions of $K$ and $B$ at the homogeneous component $\g_\gamma $.
Let $\beta \in \m^{*}$ be such that  
$$K_\gamma (X,Y)=B_\gamma (\beta_\gamma (X),Y)$$
for all $X,Y \in \g_\gamma $ and $\beta_\gamma =\beta \! \! \mid _{\g_\gamma}$. Since $B_\gamma$ is 
nondegenerate on $\g_\gamma,$ the eigenvalues of $\beta_\gamma$ are real and non zero. The eigenspaces
$\g_\gamma^1,...,\g_\gamma^p$ of $\beta_\gamma$ are pairwise orthogonal with respect to $B_\gamma$ and $K_\gamma$.
But for every $Z \in \g_e$ we have
$$K_\gamma([Z,X],Y)=K_\gamma(X,[Z,Y])=B_\gamma(\beta_\gamma(X),[Z,Y])$$
so $B_\gamma(\beta_\gamma[Z,X],Y )=B_\gamma([Z,\beta_\gamma(X)],Y )$ for every $Y \in \g_\gamma$ and 
$\beta_\gamma[Z,X]=[Z,\beta_\gamma(X)]$ that is $\beta_\gamma \circ ad \, Z=ad \, Z \circ \beta_\gamma$ 
for any $Z \in \g_e.$ This implies that $[\g_e,\g_\gamma^i] \subset \g_\gamma^i$.

\smallskip

Now we shall examinate the particular case corresponding to $\Gamma=\z$. The eigenvalues of the involutive
automorphisms $\tau_\gamma$ being real, the Lie algebra $\g$ admits a real $\Gamma$-decomposition
$\g=\sum_{\gamma \in \z}\g_\gamma.$ Then we can consider that $\g$ is a real Lie algebra.
 
Now if $i \neq j$ then
$$K_\gamma([\g_\gamma^i,\g_\gamma^j],[\g_\gamma^i,\g_\gamma^j]) \subset K([\g_\gamma^i,\g_\gamma^j],\g_e) 
 \subset (\g_\gamma^i, \g_\gamma^j)=0$$ and we have
$$[\g_\gamma^i,\g_\gamma^j]=\left\{ 0 \right\}$$
for $i\neq j$.

\medskip

\noindent {\bf Example.} In the section 4, we study the riemannian homogeneous manifold 
$SO(2l+1)/SO(r_1) \times SO(r_2)\times SO(r_3) \times SO(r_4)$. This manifold is
$(\z)$-symmetric and the Lie algebra $so \, (2l+1)$ admits a $(\z)$-grading.
By referring to the study which follows we see that
$$\g_a=A_1\oplus A_2, \ \g_b=B_1 \oplus B_2, \ \g_c=C_1\oplus C_2$$
with $[A_1,A_2]=[B_1,B_2]=[C_1,C_2]=0$ 
and we have 
$$K(A_1,A_2)=K(B_1,B_2)=K(C_1,C_2)=0.$$ 
So we have an orthogonal decomposition of each invariant space
$\g_a,\g_b,\g_c$ but the graduation is $\Gamma$-irreductible. In fact we have
$[A_1,B_1]=[A_2,B_2]=C_1, \ [A_1,B_2]=[A_2,B_1]=C_2.$ $\clubsuit $

\medskip

Let $\left\{ e,a,b,c \right\}$ be the elements of $\z$ with $a^2=b^2=c^2=e$ and $ab=c$. Each component $\g_\gamma,$ 
$\gamma \neq e$, satisfies $[\g_\gamma,\g_\gamma] \subset \g_e$ and $\g_e \oplus \g_\gamma$ is a symmetric
Lie algebra. Endowed with the inner product $\tilde{B},$ the Lie algebra $\g_e \oplus \g_\gamma$ is an 
orthogonal symmetric Lie algebra. The Killing-Cartan form is not degerate on $\g_e \oplus \g_\gamma$. 
Then $\g_e \oplus \g_\gamma$ is semi-simple. It is a direct sum of orthogonal symmetric Lie algebras of the 
following two kinds:
$$\begin{array}{ll}
i) & \g =\g' +\g' {\mbox {\rm with}}\ \g' \ {\mbox {\rm simple}} \\
ii) & \g \ {\mbox {\rm is simple.}}
\end{array}
$$
The first case has been study above and the representation is $(\z)$-irreducible.
In the second case $ad \, [\g_\gamma, \g_\gamma]$ is irreducible in $\g_\gamma$
and the representation is $(\z)$-irreducible on $\m$.

\section{Flag manifolds}
 
In this section we study riemannian properties of the oriented flag manifold
$$M=SO(2l+1)/SO(r_1) \times SO(r_2)\times SO(r_3) \times SO(r_4)$$
associated to its $\Gamma $-symmetric structures.  

\medskip
For  $\frak{g}$ classical complex simple Lie algebra of type $B_l$, it is always possible 
to endow $\g$ with a $( \z)$-grading such that
$$\g_e=so(r_1) \oplus  ... \oplus  so(r_4)$$
with $r_1+r_2+r_3+r_4=2l+1$ \cite{B.G}. The compact homogeneous space
$$M=SO(2l+1)/SO(r_1)\times SO(r_2)\times SO(r_3)\times SO(r_4)$$
is a  $(\z)$-symmetric space. We suppose $r_1\leq r_2\leq r_3\leq r_4$. 
In case  $r_1r_2 \neq 0$ and $r_3=r_4=0$ then $M$ is a symmetric space. 
The symmetric structure on the Grasmannian 
$$SO(2l+1)/SO(r_1) \times SO(r_2)$$ is well known (see \cite{K.N}).
If $r_1r_2r_3 \neq 0$, then the homogeneous space $M$ can not be symmetric. 
In what follows we shall explicitly construct on $M$ a $(\z)$-riemannian
structure. 
  Let us  consider the decomposition of a matrix of $so(2l+1)$
$$
\left(
 \begin{array}{c|c|c|c}
X_1 & A_1 & B_1 & C_1\\
\hline 
-\, ^t\! A_1 & X_2 & C_2 & B_2 \\
\hline 
-\, ^t \! B_1 & -\, ^t  C_2  & X_3 & A_2 \\
\hline 
-\, ^t  C_1 &-\, ^t \! B_2 & -\, ^t\!A_2 & X_4 \\
\end{array}
\right)
$$
with $A_1 \in \mathcal{M}(r_1,r_2), \, B_1 \in \mathcal{M}(r_1,r_3), \, C_1 \in \mathcal{M}(r_1,r_4), \,
C_2 \in \mathcal{M}(r_2,r_3), \, B_2 \in \mathcal{M}(r_2,r_4), \, A_2 \in \mathcal{M}(r_3,r_4)$ and
$X_i \in so(r_i), \ i=1,...,4.$ Let us consider the subspaces of $\g$ :
$$\g_e=
\left(
 \begin{array}{cccc}
X_1 & 0 & 0 & 0\\
0 & X_2 & 0 & 0 \\
0 & 0 & X_3 & 0 \\
0 & 0& 0 & X_4 \\
\end{array}
\right)
, \
\g_a=
\left(
 \begin{array}{cccc}
0 & A_1 & 0 & 0\\
-\, ^t\! A_1 & 0 & 0 & 0 \\
0 & 0  & 0 & A_2 \\
0&0 & -\, ^t\!A_2 & 0 \\
\end{array}
\right)
$$
$$
\g_b=
\left(
 \begin{array}{cccc}
0& 0 & B_1 & 0\\
0 & 0& 0 & B_2 \\ 
-\, ^t \! B_1 & 0  & 0&0 \\
0&-\, ^t \! B_2 & 0&0 \\
\end{array}
\right)
, \
\g_c=
\left(
 \begin{array}{cccc}
0&0&0 & C_1\\
0&0& C_2 & 0 \\
0& -\, ^t  C_2  & 0&0 \\
-\, ^t  C_1 &0&0&0 \\
\end{array}
\right).
$$
Then $\g=\g_e\oplus \g_a\oplus \g_b\oplus \g_c$ is a $(\z)$-grading of $so(2l+1)$. 
This graduation defines the $(\z)$-symmetric space
$$(SO(2l+1) ; SO(r_1)\times SO(r_2) \times SO(r_3) \times SO(r_4),{(\z)}_G).$$
Let $B$ be  a $\g_e$-invariant inner product on $\g$. By hypothesis $B(\g_\alpha ,g_\beta )=0$ as soon as $\alpha \neq \beta $ in $\z$. This shows that
$B$ is written $B=B_{\g_e}+B_{\g_a}+B_{\g_b}+B_{\g_c}$ where $B_{\g_{\alpha} }$ is an inner product on $\g_\alpha $. The restriction $B_{\g_e}$ to $\g_e$
is a biinvariant inner product. If $r_4 >2$, all the components $so(r_i)$ are simple Lie algebras and $B_{\g_e}$ is written
$$B_{\g_e}=a_1K_1+a_2K_2+a_3K_3+a_4K_4$$
where $K_i$ is the Killing-Cartan form of $so(r_i)$.  
If some components $so(r_i)$ are abelian from the index $i_0$, that is $r_i\leq 2$ for $i \geq i_0$ then
$B_{\g_e}$ is of the form $\Sigma _{j<i_0}a_jK_j \ + q$ where $q$ is a definite positive form on the abelian Lie algebra
$\oplus _{j\geq i}so(r_j)$.  Let us compute
$B_{\g_a}$. We denote by $A_1$ the subspace of $\g_a$ whose vectors are
$$
\left(
 \begin{array}{c|c|c|c}
0 & A_1 & 0 & 0\\
\hline 
-\, ^t\! A_1 & 0 & 0 & 0 \\
\hline 
 0 &0  & 0 & 0 \\
\hline 
0&0 & 0 & 0 \\
\end{array}
\right).
$$
In the same manner we define $A_2$, $B_1$, $B_2$, $C_1$ and $C_2$. For $1\leq i\leq r_1$ and $r_1+1\leq j\leq r_2$, 
let $A_{ij}$ be the corresponding elementary matrices of $A_1$
that is the $A_{ij}=(a_{rs})$ with $a_{ij}=-a_{ji}=1$ other coordinates being
equal to $0$. Similary  $X_{ij}$ denotes the elementary matrices of the  diagonal block
corresponding to $so(r_1)$, $Y_{ij}$ to $so(r_2)$, $Z_{ij}$ to $so(r_3)$ and $T_{ij}$ to $so(r_4)$. We have
$$
\left\{
\begin{array}{ll}
[X_{ij},A_{jl}]=A_{il}, & 1\leq i<j\leq r_1, \ r_1+1\leq l\leq r_1+r_2,\\
\lbrack X_{ij},A_{il} \rbrack =-A_{jl}, & 1\leq i<j\leq r_1, \ r_1+1\leq l\leq r_1+r_2,
\end{array}
\right.
$$
and
$$
\left\{
\begin{array}{ll}
[Y_{ij},A_{lj}]=A_{li}, & r_1+1\leq i<j\leq r_1+r_2, \ 1\leq l\leq r_1,\\
\lbrack Y_{ij},A_{li} \rbrack =-A_{lj}, & r_1+1\leq i<j\leq r_1+r_2, \ 1\leq l\leq r_1.
\end{array}
\right.
$$
The relation
$$B_{\g_a}([X_{rs},A_{ij}],A_{ij})=0$$
for all $X_{rs} \in so(r_1)\oplus so(r_2)$ implies
$$
\left\{
\begin{array}{lll}
B_{\g_a}(A_{ij},A_{lj})=0, & i,l\in 1,...,r_1, i\neq l, &j=r_1+1,...,r_1+r_2,\\
B_{\g_a}(A_{ij},A_{il})=0, & i=1,...,r_1, & j,l\in r_1+1,...,r_1+r_2, j \neq l.\\
\end{array}
\right.
$$ 
From the identities
$$
\left\{
\begin{array}{l}
B([X_{il},A_{lj}],A_{ij})+B(A_{lj},[X_{il},A_{ij}])=0,\\
B([Y_{ij},A_{lj}],A_{li})+B(A_{lj},[Y_{ij},A_{li}])=0,
\end{array}
\right.
$$
we obtain
$$
\left\{
\begin{array}{lll}
B_{\g_a}(A_{ij},A_{ij})= B_{\g_a}(A_{lj},A_{lj}),& i,l=1,...,r_1, &j=r_1+1,...,r_1+r_2,\\
B_{\g_a}(A_{li},A_{li})=B_{\g_a}(A_{lj},A_{lj}),  & l=1,...,r_1, & j,i=r_1+1,...,r_1+r_2.\\
\end{array}
\right.
$$ 
We deduce that all the basis vectors of $A_1$ have the same norm with respect the inner product $B$.
From the identity
$$B([X_{ij},A_{jl}],A_{js})+B(A_{jl},[X_{ij},A_{js}])=0$$
$1 \leq i<j\leq r_1, l,s \in [[ r_1+1,..., r_1 + r_2]] $, we obtain
$$B(A_{il},A_{js})+B(A_{is},A_{jl})=0.$$
Suppose that $r_1\geq 3$. There exists $r$, $1\leq r\leq r_1$  which is not equal to $i$ or $j$.
In this case we have
$$[X_{ij},A_{rs}]=0$$ and 
$$B([X_{ij},A_{jl}],A_{rs})+B(A_{jl},[X_{ij},A_{rs}])=0$$
gives
$$B(A_{il},A_{rs})=0$$
for $r \neq i$. 
This implies that the vectors $A_{ij}$ are pairwise orthogonal as soon as $r_1 >2$. It remains now to compute
$B(A_1,A_2).$ The action of $so(r_1)$ is faithful on $A_1$  and trivial on $A_2$ . Thus the $(ad_so(r_1))$-invariance of
$B_{\g_a}$ implies that $$B_{\g_a}(A_1,A_2)=0.$$ All the previous identities implies, if $r_4 >2$, that
$$B_{\g_a}=t_{A_1}\Sigma (\alpha _{ij}^1)^2+t_{A_2}\Sigma (\alpha _{ij}^2)^2,$$
where $\{\alpha _{ij}^1,\alpha _{ij}^2\}$ is the dual basis of the basis of $\g_a$ given respectively by the 
elementary matrices of $A_1$ and $A_2$ and
$t_{A_1}>0, t_{A_2}>0$. All these computations can be extended to the other components $\g_b$ and $\g_c$.
\begin{proposition}
If $r_4>2$, then all $\g_e$-invariant inner product on $\m=\g_a\oplus \g_b\oplus \g_c$ is given by
$$B=t_{A_1}\Sigma (\alpha _{ij}^1)^2+t_{A_2}\Sigma (\alpha _{ij}^2)^2+t_{B_1}\Sigma (\beta  _{ij}^1)^2+
t_{B_2}\Sigma (\beta  _{ij}^2)^2+t_{C_1}\Sigma (\gamma  _{ij}^1)^2+t_{C_2}\Sigma (\gamma  _{ij}^2)^2$$
where 
$\{\alpha _{ij}^1,\alpha _{ij}^2,\beta  _{ij}^1,\beta  _{ij}^2,\gamma  _{ij}^1,\gamma  _{ij}^2\}$ is the dual basis of 
the basis of $A_1\oplus A_2\oplus B_1\oplus B_2\oplus C_1\oplus C_2$ given by the elementary matrices and the parameters
$ t_{A_1},t_{A_2},t_{B_1},$ $t_{B_2},t_{C_1},t_{C_2}$ being nonnegative.
\end{proposition}

\medskip

\noindent It remains to examinate the particular cases corresponding to some $r_i$ equal to $2$ or $1$. This imply that
$so(r_i)$ is abelian (and not simple). 

1. If $r_1=2$ and $r_2=1$ then $r_3=r_4=1$ and the $(\z)$-grading of $so(5)$ is given by
$$so(5)=(so(2)\oplus so(1)\oplus so(1)\oplus so(1))\oplus \g_a \oplus \g_b\oplus \g_c$$ 
with dim$\g_a=3$, dim$\g_b=3$, dim$\g_c=3$ and the homogeneous space is isomorphic to 
$$SO(5)/SO(2) .$$
Every $so(2)$-invariant metric on $\m$ is of type
$$\begin{array}{ll}
B=&t_{A_1} ((\alpha _{13}^1)^2+(\alpha _{23}^1)^2)+
t_{A_2}(\alpha _{45}^2)^2+t_{B_1}( (\beta  _{14}^1)^2+(\beta _{24}^1)^2)
+t_{B_2} (\beta  _{35}^2)^2\\
&+t_{C_1}((\gamma  _{15}^1)^2+(\gamma _{25}^1)^2)+t_{C_2}(\gamma  _{34}^2)^2.
\end{array}
$$

2. If $r_1=r_2=r_3=2$ and $r_4=1$ then $\g=so(7)$. The corresponding $(\z)$-symmetric space is isomophic to
$$SO(7)/SO(2)\times SO(2)\times SO(2).$$
In this case the relation $B(A_{il},A_{rs})=0$ is not valid. We deduce that every 
$(so(2)\oplus so(2)\oplus so(2))$-invariant
inner product on $\m$ is written
$$\begin{array}{ll}
B=&t_{A_1} ( (\alpha _{13}^1)^2+(\alpha _{23}^1)^2 +(\alpha _{14}^1)^2+(\alpha _{24}^1)^2) ) +u_{A_1}
(\alpha _{13}^1\alpha _{24}^1-\alpha _{14}^1\alpha _{23}^1)\\
& +t_{A_2}((\alpha _{57}^2)^2+(\alpha _{67}^2)^2)
 +t_{B_1}( (\beta  _{15}^1)^2+(\beta _{25}^1)^2+(\beta  _{16}^1)^2+(\beta _{26}^1)^2)\\
&+u_{B_1}(\beta  _{15}^1\beta  _{26}^1-\beta  _{25}^1\beta _{26}^1)
+t_{B_2} ((\beta  _{37}^2)^2+(\beta _{ 47}^2)^2)\\
&+t_{C_1}((\gamma  _{17}^1)^2+(\gamma _{27}^1)^2)+t_{C_2}(\gamma _{36}^2)^2.
\end{array}
$$
The remaining cases correspond to $r_1=2, r_2=r_3=r_4=1$ which is treated in the example, to
$r_1=2, r_2=1, r_3=r_4=0$ and the homogeneous space is $SO(3)/SO(2)$ and it is a symmetric space and to
$r_1=r_2=r_3=1,r_4=0$ and $\g_e=\left\{ 0\right\}$. 
 So Proposition 7 and the previous results give all the
metric on flag manifolds $M$ which provide $M$ with a riemannian $(\z)$-symmetric structure. In general, for these metrics
the Levi-Civita connection is not adapted to symmetries. This connection corresponds to the canonical
torsionfree connection $\overline \nabla$ of the $(\z)$-symmetric 
homogeneous space if and only if the metric is naturally reductive
with respect to the $(\z)$-graduation. Recall that this means that
$$B([X,Y]_\m,Z)+B([X,Z]_\m,Y)=0.$$
for all $X,Y,Z \in \m$.  Applying this identity to a triple of vectors in $A_1\times B_1\times C_2$ more precisely to a
triple $(A_{r_1+1,1},B_{r_2+1,1},C_{r_1+1,r_2+1})$ we obtain that
$$t_{A_1}=t_{B_1}=t_{C_2}.$$
If we choose good triple in $A_1\times B_2\times C_2$ and $A_2\times B_2\times C_2$ we find
$$t_{A_1}=t_{B_2}=t_{C_2}$$
and
$$t_{A_2}=t_{B_2}=t_{C_2}.$$
Suppose now that the inner product corresponds to one of the  particular cases that is there is $i_0$
such that $r_{i_0}=2$. Thus in the expression of $B$ some 
double products appear. For example in the second case, $r_1=r_2=r_3=2$ and $r_4=1$. As we have
$$[B_{2,5},C_{4,5}]=-A_{2,4}$$
then
$$B(A_{1,3},[B_{2,5},C_{4,5}])+B([A_{1,3},B_{2,5}],C_{4,5}])=0$$
gives
$$B(A_{1,3},A_{2,4})=0$$
that is $u_{A_1}=0$. In the same way we find that all coefficients $u$ are equal to $0$.
\begin{proposition}
Every invariant metric $g$ on $SO(2l+1)/SO(r_1)\times SO(r_2)\times SO(r_3)\times SO(r_4)$ which is adapted 
to the $(\z)$-symmetric
structure is given by an inner product $B$ on $\m$ of type
$$B=t(\Sigma (\alpha _{ij}^1)^2+\Sigma (\alpha _{ij}^2)^2+\Sigma (\beta  _{ij}^1)^2+
\Sigma (\beta  _{ij}^2)^2+\Sigma (\gamma  _{ij}^1)^2+\Sigma (\gamma  _{ij}^2)^2)$$
with $t>0$.
\end{proposition}

\noindent {\bf Example : The homogeneous manifold $SO(5)/SO(2)\times SO(2)\times SO(1)$}

In the previous section we have described the $(\z)$-graduation of the Lie algebra $so(5)$ and we have computed the
$G$-invariant metrics which are adpated to this graduation. 
\noindent Such a metric is given by an inner product $B$ on $so(5)$ which is written
$$B=t(\alpha_1^2+ \alpha_2^2+\alpha_3^2+\alpha_4^2)+u(\alpha_1 \alpha_4-\alpha_2 \alpha_3 )+
v(\beta_1^2 +\beta_2^2)+w( \gamma_1^2+ \gamma_2^2).$$  
\begin{proposition}
\label{prop}
Every inner product on $so(5)$ for which the homogeneous components are pairwise orthogonal and which is $ad\g_e$-invariant 
is written:
$$B=q_1+t(\alpha _1^2+\alpha _2^2+\alpha _3^2+\alpha _4^2)+u(\alpha_1 \alpha_4 -\alpha_2 \alpha_3)+v(\beta_1^2+ \beta_2^2)
+w(\gamma _1^2+\gamma _2^2)$$
where $q_1$ is any inner product on $\g_e$ and $4t^2-u^2>0, \ t,v,w >0.$
This inner product gives an adapted riemannian metric on $SO(5)/SO(2)\times SO(2)$ if it is equal to
$$B=q_1+t(\alpha _1^2+\alpha _2^2+\alpha _3^2+\alpha _4^2+\beta_1^2+ \beta_2^2
+\gamma _1^2+\gamma _2^2).$$
\end{proposition}  

\noindent{\bf Remarks.} 1) If $q_1= \omega_1^2+\omega_2^2+\omega_3^2$ and $t=1,$ then $-B$ coincides 
with the Killing-Cartan form of $so(5)$. Its covariant operator $\nabla _1$ satisfies
$$2(\nabla_1) _X Y=-[X,Y].$$

2) Suppose that $g$ is the metric  $B$ on $so(5)$ corresponding to the inner product
$$B= \sum_{i=1}^3 \omega _i^2+ \sum_{i=1}^4\alpha _i^2+\sum_{i=1}^2\beta _i^2
+\sum_{i=1}^2\gamma  _i^2.$$
To simplify the notations, we shall put $E_i=A_i, \ i=1,2,3,4, \ E_5=B_1, E_6=B_2, E_7=C_1, E_8=C_2.$ 
Then the sectionnal curvatures at the identity of $\m$ are given by
$$g(R(X,Y)Y,X)_0=\frac{1}{4} B ([X,Y]_{\m},[X,Y]_{\m})+B( [X,Y]_{\g_e},[X,Y]_{\g_e})$$
and with respect to the orthonormal basis $\{E_i\}_{i=1,...,8}$ we obtain
$$
\begin{array}{l}
R_{1221}= R_{1331}=1, R_{1551}=R_{1771}=1/4 \\
R_{1441}=R_{1661}=R_{1881}=0 \\
R_{2442}=1, R_{2552}=R_{2882}=1/4 \\
R_{2332}=R_{2662}=R_{2772}=0 \\
R_{3443}=R_{3553}=R_{3883}=0 \\
R_{3663}=R_{3773}= 1/4\\
R_{4554}=R_{4774}=0 \\
R_{4664}=R_{4884}= 1/4\\
R_{5665}=1,R_{5775}=R_{5885}=1/4 \\
R_{6776}=R_{6886}=1/4 \\
R_{7887}=1.
\end{array}
$$
So the sectional curvature is positive.
 
\smallskip

3){\it On the Ambrose-Singer tensor.}

\noindent In \cite{T.V} the authors classify the homogeneous riemannian spaces using the Ambrose-Singer tensor $T$. 
The symmetric case corresponds to $T=0$. The general riemannian homogeneous spaces are classified in $8$
categories distinguished by algebraic properties of $T$. 
For the riemannian nonsymmetric space $M=SO(5)/SO(2)\times SO(2)$, this tensor corresponds to
$$T=\nabla-\overline \nabla.$$
If $\left\{ E_i \right\}_{ i=1,..,8}$ is the orthonormal basis defined above, we consider the linear map on $M$ given by
$$c_{12}(T)(X)=\sum_{i=1}^8 B_\m(T(E_i,E_j),X).$$
As $T(E_i,E_j)=-T(E_j,E_i)$, we have $c_{12}(T)(X)=0$ and $B_{\m}(T(X,Y),Z)=-B_{\m}(T(Y,X),Z)$ and
the tensor $T$ is of type ${\mathcal{T}}_3$ in the terminology of \cite{T.V}.

\medskip

4) {\it On the geodesics.} 

\noindent Following \cite{K.N}, if we set for each $X \in \m=\g_a  \oplus \g_b \oplus \g_c ,\ f_t=exp(tX) \in SO(5)$ 
and $x_t=f_t(0)\in M=SO(5)/SO(2) \times SO(2) \times SO(1)$ where $0$ is the coset $SO(2)\times SO(2) \times SO(1)$ 
 in $M$, then the curve $x_t$ is a geodesic in $M$. Conversely each geodesic starting from $0$ is of the form
$exp (tX) (0)$ for some $X \in \m$. It is not hard to see that for $E \in \m$, where $E$ stands for one of the 
$A_1,A_2,A_3,A_4,B_1,B_2,C_1$ or $C_2$, then $exp(tE)=(I_8+E^2+sintE-costE^2$ where $I_8$ is the identity
of rank $8.$ 

Two points $exp(t_1)E$ and $exp(t_2E)$ of this $2 \pi$-periodic curve falls in the same coset of $M$
if and only if $t_2-t_1=2k\pi$ for some $k \in \Z$. This shows that $f_t$ projects in a one-to-one manner 
in $M$ and its image $x_t$ is a closed geodesic (of lenght $2 \pi$).

As an example one has
$$exp(tA_1)=\left(
\begin{array}{lllll}
cost & 0 & sin t & 0 & 0 \\
0 & 1 & 0 & 0 & 0 \\
sint & 0 & cos t & 0 & 0 \\
0 & 0 & 0 & 1 & 0 \\
0 & 0 & 0 & 0 & 1 \\
\end{array}
\right).$$

\section{On lorentzian $(\z)$-symmetric structure}
It is easy to generalize the notion of riemannian $\Gamma $-symmetric homogeneous space to the notion of semi-riemannian
$\Gamma $-symmetric homogeneous space, in particular to a lorentzian metric. 
A lorentzian symmetric space $M=G/H$ is determinated by a nondegenerate $ad \h$-invariant bilinear form on $\m$ 
of signature $(1,n-1)$. In this case $M$ the Riemann curvature tensor of the Levi-Civita connection is covariant constant. 
\begin{definition}
Let $(G,H,\Gamma _G)$ a $\Gamma $-symmetric space, $g$ a semi-riemannian metric of signature $(1,n-1)$ where
$n=dim \, M$ and $B$ the corresponding 
$ad\g_e$-invariant symmetric bilinear form on $\m$ . Then  $M=G/H$ is called a $\Gamma $-symmetric lorentzian space if the homogeneous componants
of $\m$ are pairwise orthogonal with respect to $B$.
\end{definition}
Since  in the riemannian case, this doesnot imply that the riemannian connection $\nabla _g$ of $g$ coincides with $\overline \nabla$.
If $g$ satisfies this property, we will say that the connection $\nabla _g$ is adapted to the $\Gamma $-symmetric
structure. 

\noindent From the classification of $ad\g_e$-invariant form on $so(2l+1)$ given in Proposition 7,  
the $(\z)$-symmetric space
$SO(2l+1)/SO(r_1)\times ...\times SO(r_4)$ is lorentzian if and only if there exists one homogeneous component 
of $\m$ of one dimensional. For example
if we consider the $(\z)$-symmetric space
$SO(5)/SO(2)\times SO(2)\times SO(1)$ the homogeneous components are of dimension $2$ and every semi-riemannian metric is of signature
$(2p,8-2p)$ and cannot be a lorentzian metric. So $SO(5)/SO(2)\times SO(2)\times SO(1)$ can not be lorentzian. 
Nevertheless one can consider the grading of $so(5)$ given by
$$
\left(
\begin{array}{ccccc}
0 & a_1 & b_1 & b_2 & b_3 \\
-a_1 & 0 & c_1 & c_2 & c_3 \\
-b_1 & -c_1 & 0 & x_1 & x_2 \\
-b_2 & -c_2 & -x_1 & 0 & x_3 \\
-b_3 & -c_3 & -x_2 & -x_3 & 0
\end{array}
\right)
$$
where $\g_e$ is parametrized by $x_1,x_2,x_3$, $\g_a$ by $a_1$, $\g_b$ by $b_1,b_2,b_3$ and $\g_c$ by $c_1,c_2,c_3$.
Let us denote by $\{X_1,X_2,X_3,A_1,B_1,B_2,B_3,C_1,C_2,C_3\}$ the corresponding graded basis. Here $\g_e$ is
isomorphic to $so(3)\oplus so(1)\oplus so(1)$ and we obtain the $(\z)$-symmetric homogeneous space
$$SO(5)/SO(3)\times SO(1)\times SO(1)=SO(5)/SO(3).$$
Every nondegenerated symmetric bilinear form on $so(5)$ invariant by $g_e=so(3)$ is written
$$q=t(\omega _1^2+\omega _2^2+\omega _3^2)+u\alpha _1^2+v(\beta _1^2+\beta _2^2+\beta _3^2)+w(\gamma _1^2+
\gamma _2^2+\gamma _3^2)$$
where$\{\omega _i,\alpha _1,\beta _i,\gamma _i\}$ is the dual basis of the basis $\{X_i,A_1,B_i,C_i\}$.
In particular 
\begin{proposition}
The lorentzian  inner product
$$q=\omega _1^2+\omega _2^2+\omega _3^2-\alpha _1^2+\beta _1^2+\beta _2^2+\beta _3^2+\gamma _1^2+
\gamma _2^2+\gamma _3^2$$
induces a structure of lorentzian $(\z)$-symmetric structure on the nonsymmetric homogeneous space
$$SO(5)/SO(3).$$
\end{proposition}


\begin{thebibliography}{99}


\bibitem{B.G} Bahturin Y.,  Goze M., $\Gamma $-symmetric homogeneous spaces. Prepint Mulhouse (2006). 

\bibitem{Ka} Kac Victor G., {\it Infinite-dimensional Lie algebras}. Second edition. 
Cambridge University Press, Cambridge, 1985. 

\bibitem{K} Kowalski O., {\it Generalized symmetric spaces}. Lecture Notes in Mathematics, 805. 
Springer-Verlag, Berlin-New York, 1980.

\bibitem{K.N} Kobayashi S., Nomizu K., {\it Foundations of differential geometry }. Interscience tracts in pure
and applied mathematics, Number 15, Volumes I and II, Interscience publishers, 1969.  

\bibitem{L} Lutz R., {\it Sur la g\'eom\'etrie des espaces $\Gamma $-sym\'etriques}.  
C. R. Acad. Sci. Paris Sér. I Math.  293  (1981), no. 1, 55--5. 

\bibitem{L.O} Ledger A. J., Obata M., {\it Affine and Riemannian $s$-manifolds}. J. Differential Geometry  2  
1968 451--459. 

\bibitem{Rey} Reyes W., {\it A metric for a flag space}. Analysis, Geometry and Probability (Valpareiso 81), Lecture
Notes in Pure and Applied Math, 96, (1985), 173--180.

\bibitem{T.V} Tricerri F., Vanhecke L., 
{\it Homogeneous structures on Riemannian manifolds}. London Mathematical Society Lecture Note Series, {\bf 83}. 
Cambridge University Press, Cambridge, 1983.
 


\end{thebibliography}
\end{document}